\documentclass[12pt,leqno]{article}
\usepackage{amsmath, amssymb}

\newcommand{\sect}[1]{\setcounter{equation}{0}\section{#1}}

\def\epsilon{\varepsilon}

\begin{document}

\LARGE \noindent 
{\bf Holomorphic differential forms of complex manifolds 
on commutative Banach algebras and a few related problems}

\large

\vspace*{0.4em}

\hfill Hiroki Yagisita (Kyoto Sangyo University)

\vspace*{1.2em}

\normalsize

Abstract:

Let $A$ be a commutative Banach algebra. Let $M$ be a complex manifold on $A$ 
(an $A$-manifold).
Then, we define an $A$-holomorphic vector bundle 
$(\wedge^kT^*)(M)$ on $M$. For an open set $U$ of $M$, 
$\omega$ is said to be an $A$-holomorphic differential $k$-form 
on $U$, if $\omega$ is an $A$-holomorphic section 
of $(\wedge^kT^*)(M)$ on $U$. 
So, if the set of all $A$-holomorphic differential $k$-forms on $U$ 
is denoted by $\Omega_{M}^k(U)$, then $\{\Omega_{M}^k(U)\}_{U}$ 
is a sheaf of modules on the structure sheaf $O_M$ 
of the $A$-manifold $M$ 
and the cohomology group $H^l(M,\Omega_{M}^k)$ 
with the coefficient sheaf $\{\Omega_{M}^k(U)\}_{U}$ 
is an $O_M(M)$-module and therefore, in particular, 
an $A$-module. 
There is no new thing in our definition 
of a holomorphic differential form.   
However, this is necessary to get the cohomology group 
$H^l(M,\Omega_{M}^k)$ as an $A$-module. 

Furthermore, we try to define the structure sheaf 
of a manifold that is locally a continuous family of $\mathbb C$-manifolds 
(and also the one of an analytic family). 
Directing attention to a finite family of $\mathbb C$-manifolds, 
we mentioned the possibility 
that Dolbeault theorem holds for a continuous sum of $\mathbb C$-manifolds. 

Also, we state a few related problems. 
One of them is the following. 
Let $n\in \mathbb N$. 
Then, does there exist a $\mathbb C^n$-manifold $N$ 
such that for any $\mathbb C$-manifolds 
$M_1, M_2, \cdots, M_{n-1}$ and $M_n$, $N$ can not be embedded 
in the direct product $M_1\times M_2 \times \cdots \times M_{n-1} \times M_n$ 
as a $\mathbb C^n$-manifold ? 
So, we propose something that is likely to be a candidate 
for such a $\mathbb C^2$-manifold $N$.

\vfill 

\noindent 

\vfill 


%

\newpage

\normalsize

\sect{Definitions and problems} 
\ \ \ \ \ 
From $\S 2$ of [8], we follow some terms (e.g., commutative Banach algebra, 
Banach $A$-module, etc.). 

Let $A$ be a commutative Banach algebra and be fixed. 

\vspace*{0.8em}

\noindent 
{\bf Definition 1} (Topological module) : 

$X$ is said to be a topological $A$-module, 
if $X$ is an $A$-module, it is a topological $\mathbb C$-linear space, 
$$(c1_A)u = cu \ \ \ \ \ \ (c\in {\mathbb C}, \, u\in X)$$
holds and the map 
$$(a,u) \in A\times X \ \ \ \mapsto \ \ \ au \in X$$ 
is continuous. 

\noindent 
{\sf Remark} : \ 

If one is a Banach $A$-module, 
then it is a topological $A$-module. 
\hfill ---

\noindent 
{\bf Definition 2} (Linear mapping) : 

Let $X$ and $Y$ be $A$-modules. 
A mapping $F:X\rightarrow Y$ is said to be $A$-linear, 
if it satisfies 
$$F(u+v) \, = \, F(u)+F(v) \ \ \ \ \ \ (\, u, v \, \in \, X \, ),$$
$$F(fu) \, = \, fF(u) \ \ \ \ \ \ ( \, f\in A, \, \, \, u\in X \, ).$$

\noindent 
{\bf Definition 3} (Continuous multilinear mapping) : 

Let $X_1, X_2, \cdots, X_k$ and $Y$ 
be topological $A$-modules.   
A map $f$ from $X_1\times X_2\times\cdots\times X_k$ to $Y$ 
is said to be $(A,k)$-linear, 
if $f$ is $A$-linear with respect to each variable $x_i\in X_i$. 
Let $L_A \, ( \, X_1, X_2, \cdots, X_k \, ; \, Y \, )$ 
denote the set of all continuous $(A,k)$-linear mappings 
from $X_1\times X_2\times\cdots\times X_k$ to $Y$. 
$L_A \, ( \, X_1, X_2, \cdots, X_k \, ; \, Y \, )$ 
is an $A$-module. 

\noindent 
{\sf Remark} : \ 
If one is $(A,k)$-linear, 
then it is $({\mathbb C},k)$-linear. 
\hfill ---

\noindent 
{\bf Definition 4} (Norm of a multilinear mapping) : 

Let $X_1, X_2, \cdots, X_k$ and $Y$ 
be Banach $A$-modules.   
For an $(A,k)$-linear mapping $f$ 
from $X_1\times X_2\times\cdots\times X_k$ to $Y$, 
let $$\|f\| \ := \ \sup \, \{ \, \|f(x_1, x_2, \cdots, x_k)\|_Y \, |$$ 
$$\|x_1\|_{X_1}=\|x_2\|_{X_2}=\cdots=\|x_k\|_{X_k}=1 \, \}.$$ 
\hfill ---

\noindent 
{\bf Lemma 5} : 

Let $X_1, X_2, \cdots, X_k$ and $Y$ 
be Banach $A$-modules.   
Then, for any $(A,k)$-linear mapping $f$ 
from $X_1\times X_2\times\cdots\times X_k$ to $Y$,  
$f$ is continuous if and only if $\|f\|<+\infty$ holds. 
$L_A \, ( \, X_1, X_2, \cdots, X_k \, ; \, Y \, )$ 
is a Banach $A$-module.

\noindent 
{\sf Proof} : \ 
It is easy. 
\hfill 
$\blacksquare$

\noindent 
{\bf Definition 6} (Continuous antisymmetric form) : 

Let $X$ be a topological $A$-module. 
A mapping $f$ from $X^k$ to $A$ 
is said to be an $(A,k)$-linear form of $X$, 
if $f$ is $(A,k)$-linear. 
An $(A,k)$-linear form $f$ of $X$ 
is said to be antisymmetric, 
if for any permutation $\sigma$, 
$f(x_{\sigma(1)},x_{\sigma(2)},\cdots,x_{\sigma(k)}) 
\, = \, {\rm sgn}(\sigma) \cdot f(x_1,x_2,\cdots,x_k)$ 
holds. Let $A_A^k(X)$ denote the set 
of all continuous antisymmetric $(A,k)$-linear forms of $X$. 
$A_A^k(X)$ is an $A$-submodule of $L_A(X,X,\cdots,X;A)$.

\noindent 
{\bf Lemma 7} : 

Let $X$ be a Banach $A$-module. 
Then, $A_A^k(X)$ 
is a Banach $A$-submodule of $L_A(X,X,\cdots,X;A)$. 

\noindent 
{\sf Proof} : \ 
It is easy. 
\hfill 
$\blacksquare$

\noindent 
{\bf Definition 8} (Pull back) : 

Let $X_1$ and $X_2$ be topological $A$-modules. 
For a map $F$ from $X_1$ to $X_2$ 
and a map $f$ from $X_2^k$ to $A$, 
define a map $F^*(f)$ from $X_1^k$ to $A$ by 
$$(F^*(f)) \, (x_1, x_2, \cdots, x_k)
\ := \ f \, (F(x_1), F(x_2), \cdots, F(x_k)).$$ 
If $F \, \in \, \ L_A(X_1;X_2)$ 
and $f \, \in \, A_A^k(X_2)$ hold, then 
$F^*(f) \, \in \, A_A^k(X_1)$ holds.  
\hfill ---

\noindent 
{\bf Lemma 9} : 

Let $X_1$ and $X_2$ be topological $A$-modules. 
Suppose that $F$ is a bijection from $X_1$ to $X_2$. 
Suppose that $F$ and $F^{-1}$ are continuous and $A$-linear. 
Then, $F^*_{\upharpoonright A_A^k(X_2)}$ is a bijection 
to $A_A^k(X_1)$.  
$F^*_{\upharpoonright A_A^k(X_2)}$ 
and ${F^*_{\upharpoonright A_A^k(X_2)}}^{-1}$ 
are $A$-linear. 
Further, if $X_1$ and $X_2$ are Banach $A$-modules, 
then $F^*_{\upharpoonright A_A^k(X_2)}$ 
and ${F^*_{\upharpoonright A_A^k(X_2)}}^{-1}$ 
are continuous. 

\noindent 
{\sf Proof} : \ 
It is easy. 
\hfill 
$\blacksquare$

\noindent 
{\bf Definition 10} (Banach-like module) : 

$X$ is said to be a Banach-like $A$-module, 
if $X$ is a topological $A$-module 
and there exist a Banach $A$-module $Y$ 
and a bijection $F$ from $X$ to $Y$ 
such that $F$ and $F^{-1}$ are continuous 
and $A$-linear.  
\hfill ---

\noindent 
{\bf Lemma 11} : 

Let $X$ be a Banach-like $A$-module. 
Then, $A_A^k(X)$ 
is a Banach-like $A$-module. 

\noindent 
{\sf Proof} : \ 
It follows from Lemmas 7 and 9. 
\hfill 
$\blacksquare$

\newpage

\noindent 
{\bf Definition 12} (Differentiable mapping) : 

Let $X$ and $Y$ be Banach $A$-modules. 
Let $f$ be a mapping from an open set $U$ of $X$ to $Y$. 
$f$ is said to be $A$-differentiable (on $U$), 
if $f$ is Frechet differentiable 
and for any $p\in U$, the Frechet derivative $(Df)_p$ 
is $A$-linear. 
\hfill ---

\noindent 
{\bf Definition 13} (Manifold on a commutative Banach algebra) : 

Let $M$ be a Hausdorff space. Let $S$ be a set. 
$M$ is said to be an $A$-manifold 
with the system $S$ of coordinate neighborhoods, 
if the followings hold. 

For any $\varphi \in S$, there exists a Banach $A$-module $X$ 
such that $\varphi$ is a homeomorphism from an open set of $M$
to an open set of $X$. 
For any $p\in M$, there exists $\varphi \in S$ 
such that $p$ belongs to the domain of $\varphi$. 
For any $\varphi_1, \varphi_2\in S$, the coordinate transformation 
$$\varphi_2 \, \circ \, \varphi_1^{-1} \ : \ \varphi_1(U_1\cap U_2) \ 
\rightarrow \ \varphi_2(U_1\cap U_2)$$ 
is $A$-differentiable. Here, $U_1$ and $U_2$ 
are the domains of $\varphi_1$ and $\varphi_2$, respectively. 
$\mbox{}$ \hfill ---

\vspace*{0.8em}

Let $M$ be an $A$-manifold with the system $S$ of coordinate neighborhoods 
and be fixed.  
For $\varphi\in S$, we denote the domain of $\varphi$ 
by $U_\varphi$ and the Banach $A$-module 
such that $\varphi(U_\varphi)$ is an open set of it by $X_\varphi$. 
For $p\in M$, we denote the set of all $\varphi \in S$ 
such that $p\in U_\varphi$ holds by $S_p$.

\vspace*{0.4em} 

\noindent 
{\bf Definition 14} (Tangent space) : 

Let $p\in M$. As $\dot{p}_1 \sim \dot{p}_2$ indicates that 
there exist $\varphi_1, \varphi_2 \in S_p$ such that 
$\dot{p}_2=(D(\varphi_2\circ\varphi_1^{-1}))_{\varphi_1(p)}(\dot{p}_1)$
holds, $\sim$ is an equivalence relation of 
$\cup_{\varphi\in S_p}X_\varphi$. 
Let $T_p(M)$ denote the quotient set 
$(\cup_{\varphi\in S_p}X_\varphi)/\sim$. 
The tangent space $T_p(M)$ is a Banach-like $A$-module. 
\hfill ---

\noindent 
{\bf Definition 15} (Cotangent exterior space) : 

Let $p\in M$. 
Let $(\wedge^k T^*)_p (M)$ denote 
$A_A^k(T_p(M))$. The cotangent exterior space 
$(\wedge^k T^*)_p (M)$ is a Banach-like $A$-module. 
\hfill ---

\noindent 
{\bf Definition 16} (Holomorphic mapping between complex manifolds 
on a commutative Banach algebra) : 

Let $M_1$ be an $A$-manifold with a system $S_1$ of coordinate neighborhoods. 
Let $M_2$ be an $A$-manifold with a system $S_2$ of coordinate neighborhoods. 
Let $f \, : \, M_1\rightarrow M_2$ be a continuous mapping. 
$f$ is said to be $A$-holomorphic, 
if for any $\varphi_1 \in S_1$ and $\varphi_2 \in S_2$, the mapping 
$$\varphi_2 \, \circ \, f \, \circ \, \varphi_1^{-1} 
\ : \ \varphi_1(U_1\cap f^{-1}(U_2)) \ 
\rightarrow \ X_2$$
is $A$-differentiable. 
Here, $\varphi_1$ is a mapping from $U_1$ to a Banach $A$-module $X_1$ 
and $\varphi_2$ is a mapping from $U_2$ to a Banach $A$-module $X_2$. 

\noindent 
{\sf Remark} : \ 
Let $f$ be a map from an open set of a Banach $A$-module 
to a Banach $A$-module. Then, $f$ is $A$-holomorphic 
if and only if $f$ is $A$-differentiable. 
\hfill ---

\noindent
{\bf Definition 17} (Finite direct product of Banach modules) : 

Let $X_1, X_2, \cdots, X_{k-1}$ and $X_k$ 
be Banach $A$-modules. Let 
$$\|(x_1,x_2,\cdots,x_k)\| 
\ := \ \max_l \ \|x_l\|_{X_l}$$ 
$$( \, x_1\in X_1, \, x_2\in X_2, \, \cdots, \, x_k\in X_k \, ).$$ 
The finite direct product $X_1\times X_2 \times \cdots \times X_k$ 
is a Banach $A$-module. 

\noindent 
{\bf Definition 18} (Holomorphic vector bundle on a complex manifold 
on a commutative Banach algebra) : 

$E \, := \, 
(E,M,\pi,\{(U_\lambda, X_\lambda, \varphi_\lambda)\}_{\lambda\in\Lambda})$ 
is said to be an $A$-holomorphic vector bundle (on $M$), 
if it satisfies the followings. 

$E$ and $M$ are $A$-manifolds. 
$\pi$ is an $A$-holomorphic surjection from $E$ to $M$. 
Each $U_\lambda$ is an open set of $M$. 
$M \, = \, \cup_{\lambda\in\Lambda} U_\lambda$ holds. 
Each $\varphi_\lambda$ is a map 
from $\pi^{-1}(U_\lambda)$ 
to a Banach $A$-module $X_\lambda$. For $\lambda\in\Lambda$, let 
$$\pi |_\lambda \ := \ \pi_{\upharpoonright \pi^{-1}(U_\lambda)}.$$
Each map 
$$(\varphi_\lambda,\pi |_\lambda) \ : \ 
\pi^{-1}(U_\lambda) \, \rightarrow \, X_\lambda \times U_\lambda$$
is an $A$-biholomorphic map. For $\lambda\in\Lambda$ and $p\in U_\lambda$, let 
$${\varphi_\lambda} |_p 
\ := \ {\varphi_\lambda}_{\upharpoonright \pi^{-1}(\{p\})}.$$ 
For any $\lambda_1, \lambda_2 \, \in \, \Lambda$ 
and $p \, \in \, U_{\lambda_1}\cap U_{\lambda_2}$, 
the coordinate transformation 
$${\varphi_{\lambda_2}} |_p \, \circ 
\, {\varphi_{\lambda_1}} |_p^{-1} \ : \ 
X_{\lambda_1} \, \rightarrow \, X_{\lambda_2} $$
is $A$-linear. 
\hfill ---

\noindent 
{\bf Definition 19} (Tangent bundle) : 

Let $T(M)$ denote $\cup_{p\in M} \, T_p(M)$. 
\hfill ---

\noindent 
{\bf Proposition 20} : 

The tangent bundle $T(M) \rightarrow M$ 
is an $A$-holomorphic vector bundle. 

\noindent 
{\sf Proof} : \ 
It is a corollary of Proposition 34 in Section 2. 
\hfill 
$\blacksquare$

\newpage

\noindent 
{\bf Definition 21} (Cotangent exterior bundle) : 

Let $(\wedge^k T^*)(M)$ denote 
$\cup_{p\in M} \, (\wedge^k T^*)_p (M)$. 
\hfill ---

\noindent 
{\bf Proposition 22} : 

The cotangent exterior bundle $(\wedge^k T^*)(M) \rightarrow M$ 
is an $A$-holomorphic vector bundle. 

\noindent 
{\sf Proof} : \ 
It is a corollary of Proposition 38 in Section 2. 
\hfill 
$\blacksquare$

\noindent 
{\bf Theorem 23} (Probably well-known) : 

Let $A={\mathbb C}$. Let $M$ 
be an $n$-dimensional complex manifold. 
For an open set $U$ of $M$, 
the set of all holomorphic sections 
of $(\wedge^k T^*)(M)$ on $U$ is denoted by $\Omega^k_M(U)$. 
Then, the sheaf $\{\Omega^k_M(U)\}_U$ on $M$ 
is isomorphic to the sheaf of germs 
of holomorphic differential $k$-forms on $M$ 
as a sheaf of modules on the sheaf $O_M$ 
of germs of holomorphic functions on $M$. 

\noindent 
{\sf Proof} : \ 
As we define the correspondence $F^k$ by 
$$F^1(dz^i) \ : \ \ \ \ \ \ 
\sum_{j=1}^n \, \dot{z}^j \frac{\partial}{\partial z^j} 
\ \in \ T_p(M) \ \ \ \mapsto \ \ \ \dot{z}^i \ \in \ {\mathbb C}$$
$$( \, \dot{z}=(\dot{z}^1, \dot{z}^2, \cdots, \dot{z}^n) \in {\mathbb C}^n \, ),$$
$$F^k(dz^{i_1}\wedge dz^{i_2}\wedge \cdots \wedge dz^{i_k})$$
$$:= \ \sum_{\sigma \in S_k} \ {\rm sgn} (\sigma) 
\, F^1(dz^{i_{\sigma(1)}}) \otimes F^1(dz^{i_{\sigma(2)}}) 
\otimes \cdots \otimes F^1(dz^{i_{\sigma(k)}}),$$ 
it is a somewhat difficult exercise of linear algebras. 
\hfill 
$\blacksquare$

\vspace*{0.4em}

The sheaf of germs of $A$-valued $A$-holomorphic mappings on $M$ 
is denoted by $O_M$. 
From the above, we get the following definition. 

\vspace*{0.2em}

\noindent 
{\bf Definition 24} 
(Holomorphic differential form of a complex manifold 
on a commutative Banach algebra) : 

Let $U$ be an open set of $M$. 
$\omega$ is said to be an $A$-holomorphic differential $k$-form on $U$, 
if it is an $A$-holomorphic section of $(\wedge^k T^*)(M)$ on $U$. 
Let $\Omega^k_M(U)$ denote 
the set of all $A$-holomorphic differential $k$-forms on $U$. 
Then, $\{\Omega^k_M(U)\}_U$ 
is a sheaf of modules on the sheaf $O_M$ 
of germs of $A$-valued $A$-holomorphic mappings on $M$. 
The cohomology group $H^l(M,\Omega_{M}^k)$ 
with the coefficient sheaf $\{\Omega_{M}^k(U)\}_{U}$ 
is an $O_M(M)$-module and, in particular, 
an $A$-module.

\vspace*{0.2em}

\noindent
{\bf Remark} : \ 

In the coefficient sheaf $\{\Omega_{M}^k(U)\}_{U}$, 
$U$ may be limited only to all open sets of $M$ 
in an appropriate (weak) topology 
depending on the problem. 
\hfill ---

\newpage 

\noindent 
{\bf Problem 25} : 

Let $N$ be a compact continuous family 
of connected $n$-dimensional $\mathbb C$-manifolds 
on a compact Hausdorff space $X$. 
Let $\Gamma(N)$ denote the set of all continuous sections 
of $N$ on $X$. 
Then, $\Gamma(N)$ is a $C(X)$-manifold. (see [8].) 
So, if $M$ is an open set of $\Gamma(N)$, 
then the cohomology group $H^l(M,\Omega_{M}^k)$ 
is an $O_M(M)$-module and, in particular, a $C(X)$-module. 

(1) \ Let $M$ be a connected component of $\Gamma(N)$. 
Is $O_M(M)=C(X)$ ? 
Seek more specific indications 
of the sheaf $\Omega_{M}^k$ and 
the cohomology group $H^l(M,\Omega_{M}^k)$. 
When is $H^l(M,\Omega_{M}^k)$ 
a finitely generated projective $C(X)$-module ? 
Define the Euler characteristic $\chi (M,O_M)
\, = \, \sum_l \, (-1)^l[H^l(M,O_M)] \, \in \, K(X)=K(C(X))$ 
and seek its Riemann-Roch indication. 

(2) \ Let $M_1$ and $M_2$ be connected components of $\Gamma(N)$. 
Then, when are $H^l(M_1,\Omega_{M_1}^k)$ and $H^l(M_2,\Omega_{M_2}^k)$ 
isomorphic as $C(X)$-modules ?  

\noindent 
{\bf Problem 26} :

(1) \ 
Define a differential $(p,q)$-form 
and the Dolbeault operator $\overline{\partial}$ 
of an $A$-manifold. (For the case 
of an infinite-dimensional $\mathbb C$-manifold, see [3].) 

(2) \ Let $D_1, D_2, \cdots, D_{n-1}$ and $D_n$ be open disks of $A$. 
Then, is for any $l\geq1$, $H^l(D_1\times D_2\times\cdots\times D_n, 
O_{D_1\times D_2\times\cdots\times D_n})=0$ ?

(3) \ Let $U$ be a connected open set of $A^n$. 
Let $F$ be an $O_U$-module. 
Then, when is for any $l\geq1$, 
$H^l(U, F)=0$ ? For example, 
when $A=C(\{0, 1, 2, \cdots, m-1\})={\mathbb C}^m$ holds, 
define that a connected open set $U$ of ${\mathbb C}^{mn}$ 
is ${\mathbb C}^m$-Stein. Also, define 
a ${\mathbb C}^m$-coherent analytic sheaf on $U$. 

(4) \ Let $M$ be an $A$-real analytic manifold. 
Then, define the sheaf of germs of $A$-real valued 
$A$-real continuous mappings on $M$ 
and the one of $A$-real valued 
$A$-real hyperfunctions on $M$. 

\noindent
{\bf Remark} : 

Related to some of Problems 25 and 26, see Appendix 1. 
Perhaps, it may be meaningful to have $E_{k}$ 
as the sheaf of germs 
of $C^\infty$-functions on ${\mathbb R}^k$ 
in $0\rightarrow E_{n,x}\times\{0_y\}\rightarrow
E_{n,x}\times E_{m,y}\rightarrow 
\{0_x\}\times E_{m,y}\rightarrow0$. 
Apparently, a resolution that intertwined Dolbeault ones and de Rham ones   
seems to be a fine one of the structure sheaf 
of a ${\mathbb C}^m$-manifold.

In Appendix 2, 
we tried to define an analog of singular homology theory 
for a {\it continuous family of topological spaces}. 
In Appendix 3, 
we mentioned the possibility 
that Dolbeault theorem holds for a {\it continuous sum of $\mathbb C$-manifolds}. 
In Appendix 4, 
we try to define the structure sheaf 
of a manifold that is {\it locally} a continuous family of $\mathbb C$-manifolds 
(and also the one of an analytic family). 
\hfill ---

\newpage 

\noindent 
{\bf Additional Problem} :  

Consider the following $1$-dimensional $\mathbb C^2$-manifold $N$. 
It consists of four coordinate neighborhoods 
$\varphi_k:\, W_k\rightarrow \mathbb C^2 \, \, \, (k=1,2,3,4)$. 
That is, $$N=\cup_{k=1}^4 W_k.$$ Let 
$$\varphi_1(W_1):\, =\, \mathbb C\, \times \, \{\, z_2\in \mathbb C\, |\, {\rm Im} \, z_2\, >\, 0\, \},$$
$$\varphi_2(W_2):\, =\, \mathbb C\, \times \, \{\, z_2\in \mathbb C\, |\, {\rm Im} \, z_2\, <\, 0\, \},$$
$$\varphi_3(W_3):\, =\, \{\, z_1\in \mathbb C\, |\, |z_1|<1\, \}\,\times \, \mathbb C,$$
$$\varphi_4(W_4):\, =\, \{\, z_1\in \mathbb C\, |\, |z_1|<1\, \}\,\times \, \mathbb C.$$
Let $W_1\cap W_2=\emptyset$ and $W_3\cap W_4=\emptyset$. 
Let $c_1$ and $c_2$ be real numbers such that $$0<c_1<c_2<+\infty$$
holds. 
Let 
$$\varphi_3(W_1\cap W_3)
\, =\, \{\, z_1\in \mathbb C\, |\, |z_1|<1\, \}\,\times 
\, \{\, z_2\in \mathbb C\, |\, {\rm Im} \, z_2\, >\, +c_1\, \},$$
$$(\varphi_1\circ(\varphi_3^{-1}))\, (z_1,z_2)
\, =\, (\, z_1+(+2),\, z_2+(-c_1\sqrt{-1})\, ).$$
Let 
$$\varphi_4(W_1\cap W_4)
\, =\, \{\, z_1\in \mathbb C\, |\, |z_1|<1\, \}\,\times 
\, \{\, z_2\in \mathbb C\, |\, {\rm Im} \, z_2\, >\, +c_2\, \},$$
$$(\varphi_1\circ(\varphi_4^{-1}))\, (z_1,z_2)
\, =\, (\, z_1+(-2),\, z_2+(-c_2\sqrt{-1})\, ).$$
Let 
$$\varphi_3(W_2\cap W_3)
\, =\, \{\, z_1\in \mathbb C\, |\, |z_1|<1\, \}\,\times 
\, \{\, z_2\in \mathbb C\, |\, {\rm Im} \, z_2\, <\, -c_1\, \},$$
$$(\varphi_2\circ(\varphi_3^{-1}))\, (z_1,z_2)
\, =\, (\, z_1+(+2),\, z_2+(+c_1\sqrt{-1})\, ).$$
Let 
$$\varphi_4(W_2\cap W_4)
\, =\, \{\, z_1\in \mathbb C\, |\, |z_1|<1\, \}\,\times 
\, \{\, z_2\in \mathbb C\, |\, {\rm Im} \, z_2\, <\, -c_2\, \},$$
$$(\varphi_2\circ(\varphi_4^{-1}))\, (z_1,z_2)
\, =\, (\, z_1+(-2),\, z_2+(+c_2\sqrt{-1})\, ).$$

Then, answer the following question.
Do $\mathbb C$-manifolds $M_1$ and $M_2$ exist 
such that $N$ is embedded in $M_1\times M_2$ 
as a $\mathbb C^2$-manifold ? 
(If we allow that $M_1$ and $M_2$ are {\it not Hausdorff}, what about ?) 
\hfill ---

\newpage 

\sect{Proof of Propositions 20 and 22}

\noindent
{\bf Proposition 27} : 

Let $f$ be an $A$-differentiable mapping 
from an open set $U$ of a Banach $A$-module $X$ 
to a Banach $A$-module $Y$. 
Then, the mapping 
$$(\dot{z},z) \ \in \ X\times U 
\ \ \ \mapsto \ \ \ (Df)_z(\dot{z}) \ \in \ Y$$ 
is $A$-differentiable.

\noindent
{\sf Proof} : \ 
Let $F(\dot{z},z):=(Df)_z(\dot{z})$. 
Let $(\dot{z}_0,z_0) \in X\times U$. 
Then, there exists $\varepsilon>0$ such that 
$$\|z-z_0\|_X<\varepsilon 
\ \ \ \Longrightarrow \ \ \ z\in U$$ 
holds. Further, there exists $C>0$ such that 
$$\|\dot{z}_0\|_X<\frac{1}{2}\varepsilon C$$ 
holds. 
Because $f$ is $\mathbb C$-differentiable on $U$, 
$$\|z-z_0\|_X \, < \, \frac{1}{2}\varepsilon, 
\ \ \ \|\dot{z}\|_X \, < \, \frac{1}{2}\varepsilon C$$
$$\Longrightarrow$$
$$F(\dot{z},z) 
\ = \ C \ (Df)_z(\frac{1}{C}\dot{z}) \
\ = \ C \ \frac{1}{2\pi}\int_0^{2\pi}
e^{-\sqrt{-1}\theta}f(z+e^{\sqrt{-1}\theta}
\frac{1}{C}\dot{z})d\theta$$ 
holds. Hence, for any $(\dot{h},h)\in X^2$, 
$$(DF)_{(\dot{z}_0,z_0)}(\dot{h},h) 
\ = \ C \ \frac{1}{2\pi}\int_0^{2\pi} 
e^{-\sqrt{-1}\theta}(Df)_{z_0+e^{\sqrt{-1}\theta}\frac{1}{C}\dot{z}_0}
(h+e^{\sqrt{-1}\theta}\frac{1}{C}\dot{h})d\theta$$ 
holds. Because $(Df)_z$ is $A$-linear, 
for any $a\in A$ and $(\dot{h},h)\in X^2$, 
$$(DF)_{(\dot{z}_0,z_0)} \, ( \, a \, (\dot{h},h) \, )$$ 
$$= \ C \ \frac{1}{2\pi}\int_0^{2\pi} 
e^{-\sqrt{-1}\theta}(Df)_{z_0+e^{\sqrt{-1}\theta}\frac{1}{C}\dot{z}_0}
(ah+e^{\sqrt{-1}\theta}\frac{1}{C}a\dot{h})d\theta$$ 
$$= \ a \ C \ \frac{1}{2\pi}\int_0^{2\pi} 
e^{-\sqrt{-1}\theta}(Df)_{z_0+e^{\sqrt{-1}\theta}\frac{1}{C}\dot{z}_0}
(h+e^{\sqrt{-1}\theta}\frac{1}{C}\dot{h})d\theta$$ 
$$= \ a \ (DF)_{(\dot{z}_0,z_0)} \, (\dot{h},h)$$ 
holds. 
\hfill 
$\blacksquare$ 

\newpage

\noindent
{\bf Proposition 28} : 

Let $f$ be an $A$-differentiable mapping 
from an open set $U$ of a Banach $A$-module $X$ 
to a Banach $A$-module $Y$. 
Then, the mapping 
$$z \ \in \ U 
\ \ \ \mapsto \ \ \ (Df)_z \ \in \ L_A(X;Y)$$ 
is $A$-differentiable.

\noindent
{\sf Proof} : \ 
Let $G(z):=(Df)_z$. 
Let $z_0 \in U$. 
Then, there exists $\varepsilon>0$ such that 
$$\|z-z_0\|_X<\varepsilon 
\ \ \ \Longrightarrow \ \ \ z\in U$$ 
holds. Because $f$ is $\mathbb C$-differentiable on $U$, 
for any $h, \dot{z}_0 \, \in \, X$, 
$$((DG)_{z_0}(h))(\dot{z}_0)$$ 
$$= \ (1+\frac{1}{\varepsilon}\|\dot{z}_0\|_X) \, \lim_{t\rightarrow 0}
\frac{(Df)_{z_0+th}(\frac{1}{1+\frac{1}{\varepsilon}\|\dot{z}_0\|_X}\dot{z}_0)
-(Df)_{z_0}(\frac{1}{1+\frac{1}{\varepsilon}\|\dot{z}_0\|_X}\dot{z}_0)}{t}$$ 
$$= \ (1+\frac{1}{\varepsilon}\|\dot{z}_0\|_X) \frac{1}{2\pi}$$ 
$$\lim_{t\rightarrow 0} \int_0^{2\pi} e^{-\sqrt{-1}\theta} 
\frac{ 
f((z_0+th) 
+e^{\sqrt{-1}\theta}\frac{1}{1+\frac{1}{\varepsilon}\|\dot{z}_0\|_X}\dot{z}_0) 
-f(z_0 
+e^{\sqrt{-1}\theta}\frac{1}{1+\frac{1}{\varepsilon}\|\dot{z}_0\|_X}\dot{z}_0) 
}{t}
d\theta$$ 
$$= \ (1+\frac{1}{\varepsilon}\|\dot{z}_0\|_X) \, \frac{1}{2\pi}$$ 
$$\lim_{t\rightarrow 0} \int_0^{2\pi} e^{-\sqrt{-1}\theta} 
\frac{ 
f((z_0+e^{\sqrt{-1}\theta}\frac{1}{1+\frac{1}{\varepsilon}\|\dot{z}_0\|_X}\dot{z}_0) 
+th) 
-f(z_0 
+e^{\sqrt{-1}\theta}\frac{1}{1+\frac{1}{\varepsilon}\|\dot{z}_0\|_X}\dot{z}_0) 
}{t}
d\theta$$ 
$$= \ (1+\frac{1}{\varepsilon}\|\dot{z}_0\|_X) \, \frac{1}{2\pi} \, 
\int_0^{2\pi} \, e^{-\sqrt{-1}\theta} \, 
(Df)_{z_0+e^{\sqrt{-1}\theta}\frac{1}{1+\frac{1}{\varepsilon}\|\dot{z}_0\|_X}\dot{z}_0} 
(h) \, d\theta$$ 
holds. Hence, because $(Df)_z$ is $A$-linear, 
for any $a\in A$ and $h, \dot{z}_0 \, \in \, X$, 
$$( \, (DG)_{z_0}(ah) \, ) \, (\dot{z}_0)$$ 
$$= \ (1+\frac{1}{\varepsilon}\|\dot{z}_0\|_X) \, \frac{1}{2\pi} \, 
\int_0^{2\pi} \, e^{-\sqrt{-1}\theta} \, 
(Df)_{z_0+e^{\sqrt{-1}\theta}\frac{1}{1+\frac{1}{\varepsilon}\|\dot{z}_0\|_X}\dot{z}_0} 
(ah) \, d\theta$$ 
$$= \ a \, (1+\frac{1}{\varepsilon}\|\dot{z}_0\|_X) \, \frac{1}{2\pi} \, 
\int_0^{2\pi} \, e^{-\sqrt{-1}\theta} \, 
(Df)_{z_0+e^{\sqrt{-1}\theta}\frac{1}{1+\frac{1}{\varepsilon}\|\dot{z}_0\|_X}\dot{z}_0} 
(h) \, d\theta$$ 
$$= \ ( \, a \, (DG)_{z_0}(h) \, ) \, (\dot{z}_0)$$ 
holds. 
\hfill 
$\blacksquare$ 

\newpage 

\noindent 
{\bf Lemma 29} : 

Let $X_1, X_2, \cdots, X_k$ and $Y$ 
be Banach $A$-modules. 
Let $f \, \in \, L_A \, ( \, X_1, X_2, \cdots, X_k \, ; \, Y \, )$. 
Then, $f$ is $A$-differentiable. 

\noindent 
{\sf Proof} : \ 
Let $x=(x_1,x_2,\cdots,x_k) \, \in \, X_1\times X_2\times \cdots\times X_k$. 
Then, the map 
$$(h_1,h_2,\cdots,h_k) \, \in \, X_1\times X_2\times \cdots\times X_k$$
$$\mapsto \ \ \ 
f(h_1, x_2, x_3, \cdots, x_{k-1}, x_k)
+f(x_1, h_2, x_3, \cdots, x_{k-1}, x_k)$$
$$+\cdots+f(x_1, x_2, x_3, \cdots, x_{k-1}, h_k) \, \in \, Y$$ 
is $A$-linear, continuous and the Frechet derivative of $f$ at $x$.  
\hfill 
$\blacksquare$

\noindent 
{\bf Lemma 30} : 

Let $X_1$ and $X_2$ 
be Banach $A$-modules. 
Then, the map
$$(F,f) \ \in \ L_A(X_1;X_2) \times A_A^k(X_2) 
\ \ \ \ \ \ \mapsto \ \ \ \ \ \  
F^*(f) \ \in \ A_A^k(X_1)$$ 
is $A$-differentiable. 

\noindent 
{\sf Proof} : \ The map 
$$(F_1,F_2,\cdots,F_k,f) \ \in \ (L_A(X_1;X_2))^k 
\times L_A(X_2,X_2,\cdots,X_2;A)$$ 
$$\mapsto \ \ \ \ \ \  
(F_1,F_2,\cdots,F_k)^*(f) \ 
\in \ L_A(X_1,X_1,\cdots,X_1;A)$$ 
defined by 
$$((F_1,F_2,\cdots,F_k)^*(f)) \, (x_1,x_2,\cdots,x_k) 
\ := \ f \, (F_1(x_1),F_2(x_2),\cdots,F_k(x_k))$$ 
is continuous and $(A,k+1)$-linear 
and so, from Lemma 29, $A$-differentiable.  
From $F^*(f)=(F,F,\cdots,F)^*(f)$, it follows. 
\hfill 
$\blacksquare$

\noindent 
{\bf Definition 31} 
(Tangent trivialization neighborhood) : 

Let $\varphi\in S$. Let 
$(\varphi_T, \pi |_\varphi) : 
\pi^{-1}(U_\varphi) \rightarrow X_\varphi\times U_\varphi$ 
denote the local trivialization coordinate system of $T(M)$ 
corresponding to $\varphi$. 
That is, for any $\dot{p}\in \pi^{-1}(U_\varphi)$, 
$$\dot{p} 
\ = \ ( \, \{(D(\psi\circ\varphi^{-1}))_{\varphi(\pi(\dot{p}))}
(\varphi_T(\dot{p}))\}_{\psi\in S_{\pi(\dot{p})}}, 
\, \pi(\dot{p}) \, )$$
holds. 
\hfill ---

\noindent 
{\bf Definition 32} 
(Tangent open base) : 

Let $B_T$ denote the set of all sets $W$ such that 
there exist $\varphi\in S$, an open set $G$ of $X_\varphi$ 
and an open set $V$ of $U_\varphi$ 
such that $W=(\varphi_T, \pi |_\varphi)^{-1}(G\times V)$ holds. 
\hfill ---

\noindent 
{\bf Lemma 33} :  

$B_T$ is an open base of $T(M)$. 
By $B_T$, $T(M)$ is a Hausdorff space. 
For any $\varphi\in S$, the map 
$$(\varphi_T, \pi |_\varphi) \, : \  
\pi^{-1}(U_\varphi) \rightarrow X_\varphi\times U_\varphi$$ 
is a homeomorphism.

\noindent 
{\sf Proof} : \ 

$1^\circ$: \ Let $\dot{p} \in T(M)$. Then, there exists 
$\varphi \in S$ such that $\pi(\dot{p}) \in U_\varphi$ holds. 
$\dot{p} \in \pi^{-1}(U_\varphi)$ 
and $(\varphi_T, \pi |_\varphi)(\dot{p}) 
\in X_\varphi\times U_\varphi$ 
hold. So, $\dot{p} \in 
(\varphi_T, \pi |_\varphi)^{-1}(X_\varphi\times U_\varphi)$ 
holds. 

$2^\circ$: \ 
Let $\varphi_1\in S$ and $\varphi_2\in S$. 
Let $G_1$ be an open set of $X_{\varphi_1}$,  
$G_2$ be an open set of $X_{\varphi_2}$, 
$V_1$ be an open set of $U_{\varphi_1}$ 
and $V_2$ be an open set of $U_{\varphi_2}$. 
Let $\dot{p} \in (\varphi_{1T}, \pi |_{\varphi_1})^{-1}(G_1\times V_1)$ 
and $\dot{p} \in (\varphi_{2T}, \pi |_{\varphi_2})^{-1}(G_2\times V_2)$. 
Then, 
$$(D(\varphi_2\circ\varphi_1^{-1}))_{\varphi_1(\pi(\dot{p}))} 
(\varphi_{1 T}(\dot{p}))=\varphi_{2 T}(\dot{p}) \in G_2$$ 
holds. Hence, by $\pi(\dot{p}) \in V_2$ and Proposition 27, 
there exist an open set $G_0$ of $X_{\varphi_1}$ 
and an open set $V_0$ of $U_{\varphi_1}\cap U_{\varphi_2}$ 
such that $(\varphi_{1 T}(\dot{p}),\pi(\dot{p})) 
\in G_0\times V_0$ and  
$$(\dot{z},q) \, \in \, G_0\times V_0 \ \ \ \Longrightarrow \ \ \ 
( \, (D(\varphi_2\circ\varphi_1^{-1}))_{\varphi_1(q)} (\dot{z}), 
\, q \, ) \, \in \, G_2\times V_2$$
hold. Then, $G_0\cap G_1$ is an open set of $X_{\varphi_1}$, 
$V_0\cap V_1$ is an open set of $U_{\varphi_1}$ and 
$$\dot{p} \ \in \ (\varphi_{1T}, \pi |_{\varphi_1})^{-1}
((G_0\cap G_1)\times (V_0\cap V_1))$$ 
$$\subset \ 
((\varphi_{1T}, \pi |_{\varphi_1})^{-1}(G_1\times V_1)) 
\, \cap \, ((\varphi_{2T}, \pi |_{\varphi_2})^{-1}(G_2\times V_2))$$ 
holds. 

$3^\circ$: \ From $1^\circ$ and $2^\circ$, $B_T$ is an open base. 
It is easy to see that $T(M)$ is Hausdorff. 

$4^\circ$: \  Let $\varphi\in S$. 
It is easy to see that $(\varphi_T, \pi |_\varphi)$ 
is a continuous bijection. 
We show that $(\varphi_{T}, \pi |_{\varphi})^{-1}$ 
is continuous. Let $W$ be an open set 
of $\pi^{-1}(U_\varphi)$. Let $\dot{p}\in W$. 
Then, there exist $\psi\in S$, 
an open set $G$ of $X_{\psi}$  
and an open set $V$ of $U_{\psi}$ such that 
$$\dot{p} \ \in 
\ (\psi_T, \pi |_\psi)^{-1}(G\times V) 
\ \subset \ W$$ 
holds. Then, 
$$(D(\psi\circ\varphi^{-1}))_{\varphi(\pi(\dot{p}))} 
(\varphi_{T}(\dot{p}))=\psi_{T}(\dot{p}) \in G$$ 
holds. Hence, by $\pi(\dot{p}) \in V$ and Proposition 27, 
there exist an open set $G_0$ of $X_{\varphi}$ 
and an open set $V_0$ of $U_{\varphi}\cap U_{\psi}$ 
such that $(\varphi_{T}(\dot{p}),\pi(\dot{p})) 
\in G_0\times V_0$ and  
$$(\dot{z},q) \, \in \, G_0\times V_0 \ \ \ \Longrightarrow \ \ \ 
( \, (D(\psi\circ\varphi^{-1}))_{\varphi(q)} (\dot{z}), 
\, q \, ) \, \in \, G\times V$$ 
hold. Then, $G_0\times V_0$ 
is an open set of $X_{\varphi}\times U_\varphi$ and 
$$\dot{p} \ \in \ (\varphi_{T}, \pi |_{\varphi})^{-1}
(G_0\times V_0) \ \subset \ 
(\psi_{T}, \pi |_{\psi})^{-1}(G\times V) \ \subset \ W$$ 
holds. Therefore, $(\varphi_{T}, \pi |_{\varphi})^{-1}$ 
is continuous. 
\hfill 
$\blacksquare$

\noindent 
{\bf Proposition 34} :  

The Hausdorff space $T(M)$ 
is an $A$-manifold 
with $\{ \, ( \, \varphi_T, \, \varphi\circ\pi|_\varphi \, ) \, \}_{\varphi\in S}$
as the system of coordinate neighborhoods. 
The $A$-manifold $T(M)$ is an $A$-holomorphic vector bundle on $M$ 
with $\{ \, ( \, \varphi_T, \, \pi|_\varphi \, ) \, \}_{\varphi\in S}$ 
as the system of local trivialization coordinate neighborhoods.

\noindent 
{\sf Proof} : \  
It follows from Proposition 27 and Lemma 33. 
\hfill 
$\blacksquare$

\noindent 
{\bf Definition 35} 
(Cotangent exterior trivialization neighborhood) : 

Let $\varphi\in S$. Let 
$(\varphi_{\wedge^kT^*}, \pi |_\varphi) : 
\pi^{-1}(U_\varphi) \rightarrow A_A^k(X_\varphi)\times U_\varphi$ 
denote the local trivialization coordinate system 
of $(\wedge^kT^*)(M)$ corresponding to $\varphi$. 
That is, for any $f \, \in \, \pi^{-1}(U_\varphi)$ 
and $\dot{p}_1,\dot{p}_2,\cdots,\dot{p}_k \, \in \, T_{\pi(f)}(M)$, 
$$f \, (\dot{p}_1,\dot{p}_2,\cdots,\dot{p}_k)$$
$$= \ (\varphi_{\wedge^kT^*}(f)) \, ((\varphi_T|_{\pi(f)})(\dot{p}_1),
(\varphi_T|_{\pi(f)})(\dot{p}_2),\cdots,(\varphi_T|_{\pi(f)})(\dot{p}_k))$$
holds. Here, let $\varphi_T|_{p} 
\, := \, {\varphi_T}_{\upharpoonright T_{p}(M)}$ 
for $p\in U_\varphi$. 
\hfill ---

\noindent 
{\bf Definition 36} (Cotangent exterior open base) : 

Let $B_{\wedge^kT^*}$ denote the set of all sets $W$ such that 
there exist $\varphi\in S$, an open set $G$ of $A_A^k(X_\varphi)$ 
and an open set $V$ of $U_\varphi$ 
such that $W=(\varphi_{\wedge^kT^*}, \pi |_\varphi)^{-1}(G\times V)$ holds. 
\hfill ---

\noindent 
{\bf Lemma 37} : 

$B_{\wedge^kT^*}$ is an open base of $(\wedge^kT^*)(M)$. 
By $B_{\wedge^kT^*}$, $(\wedge^kT^*)(M)$ is a Hausdorff space. 
For any $\varphi\in S$, the map 
$$(\varphi_{\wedge^kT^*}, \pi |_\varphi) \, : \ 
\pi^{-1}(U_\varphi) \rightarrow A_A^k(X_\varphi)\times U_\varphi$$ 
is a homeomorphism. 

\noindent 
{\sf Proof} : \ 

$1^\circ$: \ Let $f \in (\wedge^kT^*)(M)$. Then, there exists 
$\varphi \in S$ such that $\pi(f) \in U_\varphi$ holds. 
$f \in \pi^{-1}(U_\varphi)$ 
and $(\varphi_{\wedge^kT^*}, \pi |_\varphi)(f) 
\in A_A^k(X_\varphi)\times U_\varphi$ 
hold. So, $f \in 
(\varphi_{\wedge^kT^*}, \pi |_\varphi)^{-1}
(A_A^k(X_\varphi)\times U_\varphi)$ 
holds. 

$2^\circ$: \ 
Let $\varphi_1\in S$ and $\varphi_2\in S$. 
Let $G_1$ be an open set of $A_A^k(X_{\varphi_1})$,  
$G_2$ be an open set of $A_A^k(X_{\varphi_2})$, 
$V_1$ be an open set of $U_{\varphi_1}$ 
and $V_2$ be an open set of $U_{\varphi_2}$. 
Let $f \in (\varphi_{1\wedge^kT^*}, \pi |_{\varphi_1})^{-1}(G_1\times V_1)$ 
and $f \in (\varphi_{2\wedge^kT^*}, \pi |_{\varphi_2})^{-1}(G_2\times V_2)$. 
Then, 
$$((D(\varphi_1\circ\varphi_2^{-1}))_{\varphi_2(\pi(f))})^* 
\, (\varphi_{1 \wedge^kT^*}(f)) 
\, = \, \varphi_{2 \wedge^kT^*}(f) 
\, \in \, G_2$$ 
holds. Hence, by $\pi(f) \in V_2$, Proposition 28 and Lemma 30, 
there exist an open set $G_0$ of $A_A^k(X_{\varphi_1})$ 
and an open set $V_0$ of $U_{\varphi_1}\cap U_{\varphi_2}$ 
such that $(\varphi_{1 \wedge^kT^*}(f),\pi(f)) 
\in G_0\times V_0$ and  
$$(h,q) \, \in \, G_0\times V_0 \ \ \ \Longrightarrow \ \ \ 
( \, ((D(\varphi_1\circ\varphi_2^{-1}))_{\varphi_2(q)})^* \, (h), 
\, q \, ) \, \in \, G_2\times V_2$$ 
hold. Then, $G_0\cap G_1$ is an open set of $A_A^k(X_{\varphi_1})$, 
$V_0\cap V_1$ is an open set of $U_{\varphi_1}$ and 
$$f \ \in \ (\varphi_{1\wedge^kT^*}, \pi |_{\varphi_1})^{-1}
((G_0\cap G_1)\times (V_0\cap V_1))$$ 
$$\subset \ 
((\varphi_{1\wedge^kT^*}, \pi |_{\varphi_1})^{-1}(G_1\times V_1)) 
\, \cap \, ((\varphi_{2\wedge^kT^*}, \pi |_{\varphi_2})^{-1}(G_2\times V_2))$$ 
holds. 

$3^\circ$: \ From $1^\circ$ and $2^\circ$, $B_{\wedge^kT^*}$ is an open base. 
It is easy to see that $(\wedge^kT^*)(M)$ is Hausdorff. 

$4^\circ$: \  Let $\varphi\in S$. 
It is easy to see that $(\varphi_{\wedge^kT^*}, \pi |_\varphi)$ 
is a continuous bijection. 
We show that $(\varphi_{\wedge^kT^*}, \pi |_{\varphi})^{-1}$ 
is continuous. Let $W$ be an open set 
of $\pi^{-1}(U_\varphi)$. Let $f\in W$. 
Then, there exist $\psi\in S$, 
an open set $G$ of $A_A^k(X_{\psi})$  
and an open set $V$ of $U_{\psi}$ such that 
$$f \ \in 
\ (\psi_{\wedge^kT^*}, \pi |_\psi)^{-1}(G\times V) 
\ \subset \ W$$ 
holds. Then, 
$$((D(\varphi\circ\psi^{-1}))_{\psi(\pi(f))})^* 
\, (\varphi_{\wedge^kT^*}(f)) 
\, = \, \psi_{\wedge^kT^*}(f) 
\, \in \, G$$ 
holds. Hence, by $\pi(f) \in V$, Proposition 28 and Lemma 30, 
there exist an open set $G_0$ of $A_A^k(X_{\varphi})$ 
and an open set $V_0$ of $U_{\varphi}\cap U_{\psi}$ 
such that $(\varphi_{\wedge^kT^*}(f),\pi(f)) 
\in G_0\times V_0$ and  
$$(h,q) \, \in \, G_0\times V_0 \ \ \ \Longrightarrow \ \ \ 
( \, ((D(\varphi\circ\psi^{-1}))_{\psi(q)})^* \, (h), 
\, q \, ) \, \in \, G\times V$$ 
hold. Then, $G_0\times V_0$ 
is an open set of $A_A^k(X_{\varphi})\times U_\varphi$ and 
$$f \ \in \ (\varphi_{\wedge^kT^*}, \pi |_{\varphi})^{-1}
(G_0\times V_0) \ \subset \ 
(\psi_{\wedge^kT^*}, \pi |_{\psi})^{-1}(G\times V) \ \subset \ W$$ 
holds. Therefore, $(\varphi_{\wedge^kT^*}, \pi |_{\varphi})^{-1}$ 
is continuous. 
\hfill 
$\blacksquare$

\newpage

\noindent 
{\bf Proposition 38} :  

The Hausdorff space $(\wedge^kT^*)(M)$ 
is an $A$-manifold 
with $\{ \, ( \, \varphi_{\wedge^kT^*}, \, 
\varphi\circ\pi|_\varphi \, ) \, \}_{\varphi\in S}$ 
as the system of coordinate neighborhoods. 
The $A$-manifold $(\wedge^kT^*)(M)$ is an $A$-holomorphic vector bundle on $M$ 
with $\{ \, ( \, \varphi_{\wedge^kT^*}, \, \pi|_\varphi \, ) \, \}_{\varphi\in S}$ 
as the system of local trivialization coordinate neighborhoods. 

\noindent 
{\sf Proof} : \  
It follows from Proposition 28 and Lemmas 30, 37. 
\hfill 
$\blacksquare$ 


\vfill 

Acknowledgment: 

This work was supported by JSPS KAKENHI Grant Number JP16K05245. 

\noindent 
ArXiv does not seem to accept frequent revisions. 
The revised version may have been put in 
``\verb|https://www.researchgate.net/profile/Hiroki_Yagisita|''.

\newpage

{\bf Appendix}

\vspace*{1.6em}

\noindent 
{\sf 1} : \  
Let $X$ be a compact Hausdorff space. 
Let $U$ be a convex open set of ${\mathbb R}^n$. 
Let $F$ be a $C^1$-map from $C(X;U)$ to $C(X;{\mathbb R})$. 
Suppose that for any $u \in C(X;U)$, the Frechet derivative 
$F^\prime (u)$ of $F$ at $u$ is $C(X;{\mathbb R})$-linear. 
Then, for any $u_0, u_1 \in C(X;U)$ and $x \in X$,  
$u_0(x)=u_1(x)$ implies $(F(u_0))(x)=(F(u_1))(x)$. 
So, there exists a function $f$ from $X\times U$ to $\mathbb R$ 
such that for any $u \in C(X;U)$ and $x \in X$, 
$(F(u))(x)=f(x,u(x))$ holds. 

\noindent 
{\sf Proof} : \  From 
$$F(u_1)-F(u_0) \ = \ \int_0^1 \, (F^\prime((1-t)u_0+tu_1)) \, (u_1-u_0) \, dt$$
$$= \ \int_0^1 \, (F^\prime((1-t)u_0+tu_1)) \, 
(\sum_{k=1}^n (u_1^{(k)}-u_0^{(k)})e_k) \, dt$$
$$= \ \sum_{k=1}^n \ 
\int_0^1 \, (u_1^{(k)}-u_0^{(k)}) \, (F^\prime((1-t)u_0+tu_1))(e_k) \, dt,$$
$$(F(u_1))(x)-(F(u_0))(x)$$
$$= \ \sum_{k=1}^n \ 
\int_0^1 \, (u_1^{(k)}(x)-u_0^{(k)}(x)) \, ((F^\prime((1-t)u_0+tu_1))(e_k))(x) \, dt$$
$$= \ \sum_{k=1}^n \ 
\int_0^1 \, 0 \, ((F^\prime((1-t)u_0+tu_1))(e_k))(x) \, dt \ = \ 0$$
holds. 
\hfill 
$\blacksquare$ 


\noindent 
{\sf 2} : \ I try to define an analog of singular homology theory 
for a continuous family of topological spaces, 
but I do not know whether it will work or not. 

Let $\pi$ be a continuous mapping 
from a topological space $M$ to one $X$. 
Let $\Delta^k$ denote the standard $k$-simplex.  
Let $T^\pi_k$ be the set of all continuous mappings 
$\sigma: \, \Delta^k \times X \rightarrow M$ 
such that for any $(x,t)\in\Delta^k \times X$, 
$$(\pi\circ\sigma)(x,t)=t$$ holds. 
Then, let $S^\pi_k$ be the free $\mathbb Z$-module 
such that $T^\pi_k$ is a base of $S^\pi_k$. 
In an appropriate class, it appears to be isomorphic 
to the {\it usual} singular homology. 
$\mbox{\, \, \,}$  \hfill ---

\newpage

\noindent 
{\sf 3} : \ 
I am considering the following in the interview, 
but I do not know whether it will work or not. 

\noindent 
[Title] \ Dolbeault theorem for a topological sum of complex structures. 

\noindent 
[Abstract] \ For an open set $U$ of ${\mathbb C}^n\times{\mathbb R}^m$, 
let $O(U)$ denote the ring 
of all $\mathbb C$-valued continuous functions $f(z,t)$ on $U$ 
such that $f(z,t)$ is holomorphic 
with respect to the several complex variables $z \in {\mathbb C}^n$. 
Then, $\{O(U)\}_U$ is a sheaf of commutative rings 
on ${\mathbb C}^n\times{\mathbb R}^m$. 
We show Dolbeault theorem 
$$H^q(U,O_U)\cong H_{\overline{\partial}_z}^q(U,O_U).$$ 
Further, we show that if $D$ is an open polydisk of ${\mathbb C}^n$, 
$T$ is an open set of ${\mathbb R}^m$ and $U=D\times T$ and $q\geq 1$ hold, 
then this cohomology is vanishing. 
So, it is shown that 
for a continuous family of additive Cousin date on an open polydisk, 
there exists a continuous family of solutions of the first problem.

\noindent 
[Comments] \ An open set $U$ of ${\mathbb C}^n\times{\mathbb R}^m$ 
is a very simple example of a continuous sum of $\mathbb C$-manifolds. 
On the other hand, an appropriate Banach $\mathbb C$-manifold 
is a continuous product. 
A projection $\pi: U\rightarrow {\mathbb R}^m$ 
is a simple example of a continuous family. 
It seems that a continuous sum and a continuous family 
have not been fully studied yet. 
A finite sum and a finite product are familiar, 
but it seems that a {\it finite family}  
is not paying much attention. 
Although a continuous sum and a continuous family 
are unexplored places, they may remain unexplored 
for a while in the future or may be stepped on in a moment. 

The $C(X)$-manifold $M$ 
corresponding to a continuous family $\pi : \, U\rightarrow X$ 
and the total space $U$ 
do not seem to reflect all of the important structures of $\pi$. 
It is desirable to construct a consistent structure as a whole 
by adding some additional information to $M$ or $U$. 
For example, it seems that the sheaf $\{M|_D\}_D$ on $X$ 
reflects most information on $\pi$. 
Here, for an open set $D$ of $X$, 
$M|_D$ is the $C_b(D)$-manifold corresponding to $\pi |_D$. 
$\mbox{\, \, \,}$  \hfill ---


\noindent 
{\sf 4} : \ 
We try to define the structure sheaf 
of a manifold that is {\it locally} a continuous family of $\mathbb C$-manifolds 
(and also the one of an analytic family). 

\noindent 
[Continuous family] \ Let $\pi: \, M \rightarrow X$ be a continuous family 
of $\mathbb C$-manifolds. Then, for an open set $U$ of $M$, let $O_{\pi}(U)$ 
be the set of all pairs $(f,g)$ of $\mathbb C$-valued continuous functions on $U$ 
such that for any $t\in X$, $f_{\upharpoonright \pi^{-1}(\{t\})}$ is locally constant 
and $g_{\upharpoonright \pi^{-1}(\{t\})}$ is $\mathbb C$-holomorphic. 

\noindent 
[Analytic family] \ Let $\pi: \, M \rightarrow X$ be an analytic family 
of $\mathbb C$-manifolds. Then, for an open set $U$ of $M$, let $O_{\pi}(U)$ 
be the set of all pairs $(f,g)$ of  $\mathbb C$-holomorphic functions on $U$ 
such that for any $t\in X$, $f_{\upharpoonright \pi^{-1}(\{t\})}$ is locally constant. 
\hfill ---


\noindent 
{\sf +0} : \ 
Let $M \rightarrow X$ be {\it an analytic family} 
of compact $\mathbb C$-manifolds. Let $U$ 
be a relatively compact open set of $X$. 
Let $\Gamma(M|U)$ denote the set of all continuous sections 
of $M$ on $U$. Then, $\Gamma(M|U)$ is a $C_b(U)$-manifold. (see [8].) 
However, perhaps, $\Gamma(M|U)$ can be considered 
as a manifold {\it with a stronger structure}. 
The following problem is related to it. 
Suppose that $X$ is Stein and $\Gamma(M)=(C_b(X))^n$ holds. 
If $\pi_1: \, M \rightarrow X$ and $\pi_2: \, M \rightarrow X$ 
are analytic families, does $\pi_1=\pi_2$ hold ? 
(Oka-Grauert principle)
\hfill ---

\noindent 
{\sf +1} : \ 
Let $X$ be a compact Hausdorff space. 
Then, define a $C(X)$-analytic set of a Banach $C(X)$-module.
Also, define a $C(X)$-algebraic set of a finitely generated projective $C(X)$-module.
\hfill ---

\noindent 
{\sf +2} : \ 
Define a locally direct product space
of complex analytic spaces
and its underlying complex analytic space.
Define a locally direct product space
of complex algebraic varieties (or, schemes etc.)
and its underlying complex algebraic variety (or, scheme etc.).
However, the rudimentary general theory of a ringed space
may bring a locally direct product space.
On the other hand, an underlying space may be more difficult.
\hfill ---

\vfill


\vspace*{3.2em}

{\bf References}

\vspace*{0.8em}

[1] S. Araki, Topological K-theory (Japanese), {\it Sugaku}, 22 (1970), 60-76. 

[2] L. Hormander, $L^2$ estimates and existence theorems 
for the $\overline{\partial}$ operator, {\it Acta Math.}, 113 (1965), 89-152. 

[3] L. Lempert, The Dolbeault complex in infinite dimensions, {\it J. Amer. Math. Soc.}, 11 (1998), 485-520. 

[4] A. Mallios and E. E. Rosinger, 
Space-time foam dense singularities and de Rham cohomology, 
{\it Acta Appl. Math.}, 67 (2001), 59-89. 

[5] T. Ohsawa and K. Takegoshi, On the extension of $L^2$ holomorphic functions, 
{\it Math. Z.}, 195 (1987), 197-204. 

[6] H. Ozeki, Vector bundles and projective modules (Japanese), {\it Sugaku}, 18 (1967), 223-233. 

[7] M. H. Papatriantafillou, 
Partitions of unity on $A$-manifolds, 
{\it Internat. J. Math.}, 9 (1998), 877-883. 

[8] H. Yagisita, Finite-dimensional complex manifolds on commutative Banach algebras and continuous families of compact complex manifolds, {\it arXiv.org}. 


\end{document}